\documentclass[winedt,yap]{iitparc}
\usepackage{amssymb,amsmath}
\usepackage{amsfonts}
\usepackage[mathscr]{eucal}

\def\x#1{} 
\def\cdc{,\ldots,}
\def\la{\lambda}
\def\ind{\mathop{\rm ind}\nolimits}
\def\adj{\mathop{\rm adj}\nolimits}
\def\rank{\mathop{\rm rank}\nolimits}
\def\CC{{\mathbb C}}
\def\nn{{n\times n}}
\begin{document}

\frontmatter

\JournalISSNCode{0005-1179}
\OrigYearOfIssue{2002}%
\OrigIssueNo{10}%
\IssuePrice{27.00}%
\TransYearOfIssue{2002}%
\TransCopyrightYear{2002}%
\OrigCopyrightYear{2002}%
\TransVolumeNo{63}%
\TransIssueNo{10}%

\mainmatter

\setcounter{page}{1537}

\Rubrika{DETERMINATE SYSTEMS}

\newcommand\propro[2]{\par\PPR{#1}{\rm #2}\medskip}
\newcommand\prothe[2]{\par\PTH{#1}{\rm #2}\qed\medskip}
\newcommand\prolem[2]{\medskip\par\PLE{#1}{\rm #2}\qed\medskip}
\newcommand\procor[2]{\medskip\par\PCO{#1}{\rm #2}\qed\medskip}

\title{On Determining the Eigenprojection\\ 
         and Components of a Matrix%
\thanks{This work was supported by the Russian Foundation for Basic Research, project no.~02--01--00614-a.}}

\titlerunning{On Determining the Eigenprojection and Components}

\author{R. P. Agaev and P. Yu. Chebotarev}
\authorrunning{Agaev, Chebotarev}
\OrigCopyrightedAuthors{Agaev, Chebotarev}
\institute{Trapeznikov Institute of Control Sciences, Russian Academy of Sciences, Moscow, Russia}
\received{Received June 6, 2002 {\bf[a revised version of January 21, 2011]}.}
\OrigPages{3--12}

\maketitle

\begin{abstract}
Matrix theory and its applications make wide use of the eigenprojections of square matrices. The present\x{} paper demonstrates\x{} that the eigenprojection of a matrix $A$ can be calculated with the use of any annihilating polynomial of\x{} $A^u,$ where $u\ge\ind A.$ This enables one to find\x{} the components and the minimal\x{} polynomial of $A,$ as well as the Drazin inverse~$A^D.$
\end{abstract}

\section{introduction}\label{s-intro}\setcounter{footnote}{1}

\subsection{On Applications of Drazin Inverse, Eigenprojection, and Matrix Components}\vspace*{3mm}

The Drazin inverse, eigenprojections, and components of square matrices, as well as their minimal\x{} polynomials, are widely used in mathematics, mechanics, engineering, economics, and other sciences. Let us mention some of these applications. Dynamics of multibody systems is described by systems of differential-algebraic equations. Their linearization results in ``implicit'' singular systems of ordinary differential equations. Their analysis relies on determining the Drazin inverse of\x{} the matrices of coefficients~\cite{CampbellMeyerRose76}. This approach was used in \cite{SimeonFuhrerRentrop93} to analyze the stabilization of a model of truck motion. The eigenprojections that correspond to the eigenvalues of a matrix are the main components of this matrix.\x{} From the computational point of view, determination of the eigenprojection corresponding to the eigenvalue $0$ amounts to determining the Drazin inverse of the matrix\x{} and it can also\x{} be used to determine the Moore--Penrose inverse. The matrix components and the functions of matrices determined with their use are employed in the dynamics of rigid bodies \cite{Macmillan51}, motion stability analysis, and integration of systems of linear differential equations \cite{Gantmacher66,Lankaster78}. Applications of the eigenprojection and Drazin inverse in the theory of Markov chains \cite{Meyer75}, graph theory \cite{CheAga02a+}, cryptography \cite{HartwigLevine81}, and econometrics \cite{Chipman76} (see also~\cite{Wilkinson82}) deserve mentioning as well. 
The present paper demonstrates\x{} that for determining the eigenprojections, the Drazin inverse, and the matrix components one may use any annihilating polynomial of\x{} $A^u,$ where $A$ is a given matrix and $u\ge\ind A.$

\vspace*{2mm}\subsection{Basic Definitions}\vspace*{1mm}

Let $A\in\CC^\nn$ be a square matrix whose range and null space are denoted by ${\mathcal R}({A})$ and ${\mathcal N}({A}),$ respectively. Let $r$ be the total multiplicity of the nonzero eigenvalues of $A,$\, $v=n-r$ be the multiplicity of $0$ as the eigenvalue of $A,$ and $\nu=\ind A$ be the {\em index\/} of $A,$ that is, the smallest ${k}\in\{0,1,\ldots\}$ for which\; $\rank A^{k+1} =\rank A^{k}\,$~\cite{Ben-IsraelGreville74,CampbellMeyer79}. The index is $0$ if and only if $A$ is nonsingular. The index of a singular matrix is the index of its eigenvalue $0,$ that is, the multiplicity of zero as the root of its minimal\x{} polynomial, or the size of the greatest block with zero diagonal in its Jordan form.

The {\em eigenprojection}%
\footnote[2]{The eigenprojection is also called the {\em principal idempotent} (see, for example,~\cite{Hartwig76}).}
of a matrix $A$ {\em corresponding to the eigenvalue} $0$ \cite{Rothblum76a} or, for brevity,\x{} the {\em eigenprojection of\/} $A$ \cite{Rothblum76} is\x{} the projection, i.e., the idempotent matrix, $Z$ such that
\begin{gather*}
{\mathcal R}(Z)={\mathcal N}(A^\nu)\quad \text{and}\quad {{\mathcal N}}(Z)={\mathcal R}(A^\nu).
\end{gather*}
$Z$ can be said to be the projection {\em on ${\mathcal N}(A^\nu)$ along ${\mathcal R}(A^\nu).$} The eigenprojection is unique because an idempotent matrix is uniquely determined by its range and null space (see, for example, \cite[p.~50]{Ben-IsraelGreville74}).

The eigenprojections are used in the theory of generalized inverses and numerous applications of linear algebra, as well as for decomposing matrices into {\em components}, which allows one to determine the values of $f(A)$ for a wide class of functions $f:{\CC\to\CC}$ (see \cite{Macmillan51,Gantmacher66,Lankaster78}). In this connection, we also note that for any finite homogeneous Markov chain with the transition\x{} matrix $P,$ the limiting matrix of mean probabilities
\begin{gather*}
P^{\infty}=\lim_{k\to\infty}\frac{1}{k}\,\sum_{t=0}^{k-1} P^t
\end{gather*}
is the eigenprojection of the matrix $I-P$~\cite{Meyer75,Rothblum76a}. Another example: for any weighted digraph $\Gamma$ with the Laplacian matrix $L,$ the matrix $\Gamma$ of maximum converging forests is the eigenprojection of the matrix $L$~\cite{CheAga02a+}. The present paper is devoted to some methods of determining the eigenprojection of a matrix, as well as its components and the minimal\x{} polynomial.

\vspace*{2mm}\subsection{Preliminaries}\vspace*{1mm}

From the fact that $\nu$ is the order of the greatest Jordan block with zero diagonal in the Jordan form of $A,$ it follows that
\begin{gather}
\label{nu<n} \nu\le v.
\end{gather}
For the powers of a Jordan block $J$ with zero diagonal we get
\begin{gather*}
\rank J^k=\max\{l-k,0\},
\end{gather*}
where $l$ is the order of $J.$ Therefore,
\begin{gather}
\label{e-rankA^u=r}\rank A^i>\rank A^{i+1}\quad\mbox{for}\quad 0\le i<\nu\quad\mbox{ and}\quad\rank A^i=r\quad\mbox{for}\quad i\ge\nu\quad \mbox{and}
\\
\label{e-indA^k=1}\ind A^k=\min\,\{t\in{\mathbb Z}\mid kt\ge\nu\}=\lceil\nu/k\rceil,
\end{gather}
where $\lceil x\rceil$ is the smallest integer not smaller than $x.$

Since ${\mathcal R}(X)$ and ${\mathcal N}(X)$ are subspaces of $\CC^n$ for any matrix $X\in\CC^\nn$ and since
\begin{gather*}
\dim{\mathcal R}(X)=n-\dim{\mathcal N}(X)=\rank X,
\end{gather*}
we obtain by using \eqref{e-rankA^u=r} that
\begin{gather}
\label{NN(A)} {\mathcal N}(A)\subset\ldots\subset{\mathcal N}(A^\nu)={\mathcal N}(A^{\nu+1})=\ldots,
\\
\label{RR(A)} {\mathcal R}(A)\supset\ldots\supset{\mathcal R}(A^\nu)={\mathcal R}(A^{\nu+1})=\ldots,
\end{gather}
where the strict inclusions refer to the case of $\nu>1.$ These relations and \eqref{e-indA^k=1} give rise to the following lemma.

\begin{lemma}
\label{l-prs} 
The eigenprojections of the matrices $A,A^2,\ldots A^\nu,\ldots$ are the same.
\end{lemma}

\vspace*{2mm}\subsection{Plan of the Paper}\vspace*{1mm}

The paper makes use of the following approach. The eigenprojection of a matrix with index\footnote{In this case, the eigenprojection is the zero matrix.\x{}} $0$ or $1$ is the easiest to determine. Therefore, knowing an upper bound $u$ for $\,\ind A,$\x{} one can seek, by virtue of \eqref{e-indA^k=1}, the eigenprojection of $A^u$ (which has index~$1$) instead of the eigenprojection of~$A.$ Indeed, according to Lemma~\ref{l-prs}, the eigenprojections of $A$ and $A^u$ coincide. In the next section, we use this approach and several known\x{} characterizations of the eigenprojection to obtain simple expressions for the eigenprojection of $A^u$ and, consequently, of~$A.$ In Section~\ref{s-annih}, we propose a method for determining the eigenprojection of $A$ (and the Drazin inverse $A^D$) by means of any nonzero annihilating polynomial of\x{}~$A^u.$ Section~\ref{s-co_mi} discusses calculation of the components and the minimal\x{} polynomial of a square matrix. Explicit expressions for\x{} the eigenprojection and the components of a matrix $A$ in terms of its eigenvalues are given in Section~\ref{s-ep_eva}.

\vspace*{3mm}\section{characterizations of the eigenprojection}\label{s-chara}\vspace*{2mm}

Below we present a number of conditions that are equivalent to the fact that an idempotent matrix $Z\in\CC^\nn$ is the eigenprojection of a matrix $A\in\CC^\nn$ having index~$\nu.$ This list refers also to the publications where the corresponding\x{} conditions are presented.
\begin{itemize}
\item[(a)] ${\mathcal R}(Z) ={\mathcal N}(A^\nu)$ and ${\mathcal R}(Z^*)={\mathcal N}((A^*)^\nu),$ where $X^*$ is the Hermitian adjoint of $X$~\cite{Rothblum76a}.

\item[(b)] $A^\nu Z=ZA^\nu=0$ and $\rank A^\nu+\rank Z=n$~\cite{Wei96,Zhang01}.

\item[(c)] $AZ=ZA$ and $A+\alpha Z$ is nonsingular for all $\alpha\ne0$ \cite{KolihaStraskraba99}.

\item[(d)] $AZ=ZA,$\, $A+\alpha Z$ is nonsingular for some $\alpha\ne0,$ and $AZ$ is nilpotent \cite{KolihaStraskraba99}.

\item[(e)] $AZ=ZA,$\, $AZ$ is nilpotent, and $AU=I-Z=VA$ for some $U,V \in\CC^\nn$~\cite{Harte84}.

\item[(f)] $ZC=CZ$ for any $C$ commuting with $A,$\, $AZ$ is nilpotent, and $(\det A =0\Rightarrow Z\ne0)$~\cite{Koliha01}.
\end{itemize}

We now turn\x{} to some more constructive characterizations of the eigenprojection. In the items to follow, $u$ is any integer  not less than\x{!!} $\nu,$ and we use Lemma~\ref{l-prs} according to which $A$ and $A^u$ have the same eigenprojection. Also, the following notation is used: $A^D$ is the {\em Drazin inverse\/} of\x{}~$A,$\, $(A^u)^{\scriptscriptstyle\#}$ is the {\em group inverse\/} of\x{} $A^u,$
\begin{gather}
\label{e-charA} \psi(\la)=\sum_{j=v}^{n}a_{n-j}\la^j
\end{gather}
and
\begin{gather}
\label{e-adjA} B(\la)=\adj(\la I-A)=\sum_{j=0}^{n-1}A_{n-1-j}\la^j
\end{gather}
are, respectively, the {\em characteristic polynomial\/} of\x{} $A$ and the {\em adjoint matrix\/} of\x{} $\la I-A;$ $\psi_u(\la)$ and $\widehat\psi_u(\la)$ are, respectively, the {\em characteristic\/} and {\em minimal\x{} polynomials\/} of\x{} $A^u;$ $d_u(\la)$ is their quotient:
\begin{gather}
\label{e-phi(la)/d} \psi_u(\la)=d_u(\la)\,\widehat\psi_u(\la);
\end{gather}
$B_u(\la)=\adj(\la I-A^u)$ is the {\em adjoint matrix\/} of\x{} $\la I\!-\!A^u,\x{}$ and $C_u(\la)$ is the {\em reduced adjoint matrix\/} \cite{Gantmacher66} of\x{} $\la I-A^u,$ that is, the
matrix polynomial
\begin{gather}
\label{e-C(la)} C_u(\la)=\frac{B_u(\la)}{d_u(\la)},
\end{gather}
which can also be obtained by dividing the matrix polynomial $\widehat\psi_u(\la)I$\x{??} by the binomial $\la I-A^u.$ Some constructive characterizations of the eigenprojection are as follows.
\begin{itemize}
\item[(g)] $Z=I-AA^D=I-A^u(A^u)^{\scriptscriptstyle\#}$
\cite{Rothblum76,Hartwig76}.

\item[(h)] $Z=I-(I-A_{n-v}/a_{n-v})^u$ \cite{Hartwig76}.%
\footnote[3]{The coefficients $a_i$ of the characteristic polynomial \eqref{e-charA} and the matrix coefficients $A_i$ of the adjoint matrix \eqref{e-adjA} can be calculated concurrently by Faddeev's algorithm \cite{FaddeevFaddeeva63R,Gantmacher66}.}\x{}
%

\item[(i)] ${Z={C_u(0)\big/\rho_u(0)}},$ where $\rho_u(\la)$ is the polynomial satisfying\x{} the equality $\la\,\rho_u(\la)= \widehat\psi_u(\la)$~\cite{Gantmacher66}.

\item[(j)] ${Z=\lim\limits_{\la\to0}\la B_u(\la)\big/\psi_u (\la)},$ where $\la$ is such that\x{} $\psi_u(\la)\ne0.$ This follows from \eqref{e-phi(la)/d}, \eqref{e-C(la)}, and~(i).

\item[(k)] $Z=\lim\limits_{\left|\tau\right|\to\infty}(I+\tau A^u)^{-1}.$ 
This follows from item\x{} (j) or Theorem~3.1 in\x{}~\cite{Meyer74}.

\item[(l)] $Z=X(Y^*X)^{-1}Y^*,$ where $X$ and $Y$ are the matrices whose columns make up bases in ${\mathcal N}(A^\nu)$ and ${\mathcal N}((A^*)^\nu),$ respectively \cite{Rothblum76a,Hartwig76} ($\nu$ can be replaced by $u$ by virtue of \eqref{NN(A)} and \eqref{RR(A)}).
\end{itemize}

In the following section, we propose one more\x{another} constructive characterization of the eigenprojection. We do not touch upon the issues of computational efficiency and just note that the corresponding algorithm does not require computation of limits, the reduced adjoint matrix, the matrix coefficients of the adjoint matrix, generalized inverses, or bases of subspaces. If an arbitrary annihilating polynomial of\x{} $A^u$ or the eigenvalues of $A$ are known, then the eigenprojection of $A$ can be expressed as a scalar polynomial of~$A^u.$ If the characteristic or minimal\x{} polynomial is used as the annihilating polynomial, then the expression we propose in the next section is equivalent to the ones\x{} presented in (j) and~(i), respectively.

\vspace*{3mm}\section{calculating the eigenprojection of $A$\\ Using\x{} any annihilating polynomial of\x{} $A^u$} \label{s-annih}\vspace*{3mm}

The main result of this section is\x{} Theorem~\ref{t-eigen0} which states\x{} that the eigenprojection of\x{} $A$ can be found by direct calculation using any annihilating polynomial of\x{} $A^u$ with any $u\ge\nu=\ind A.$ One can always substitute $v$ for $u$ (see \eqref{nu<n}) and the characteristic polynomial%
\footnote[4]{It can be calculated using the Leverrier--Faddeev algorithm, the Krylov method, or other computational algorithms (see, for example, \cite{FaddeevFaddeeva63R,Gantmacher66}\x{}).}
for\x{} the annihilating polynomial, but to reduce computations, it would be advantageous to make use of the information available about $A$ for getting a more accurate upper estimate of\x{} $\nu$ and, perhaps, for obtaining an annihilating polynomial of a power lower than that of the characteristic polynomial. In this respect, the minimal\x{} polynomial of\x{} $A^\nu$ is the best choice, but its determination can be problematic itself.

The proof of Theorem~\ref{t-eigen0} relies\x{} on the above condition~(f) which is necessary and sufficient for a matrix $Z$ to be the eigenprojection of~$A.$ We note, however, that employing\x{} the definition of the eigenprojection would not complicate the proof considerably.

Let, as before, $u\ge\nu$ (for example, $u=v$ or $u=n$) and let $\varphi(\la)$ be an arbitrary nonzero annihilating polynomial of\x{} $A^u$: $\varphi(\la)\equiv0$ is not true and $\varphi(A^u)=0.$ Let
\begin{gather}
\label{eigpr2} \varphi(\la)=\la^t(\la^q+p_1\la^{q-1}+\ldots +p_q),\quad
\end{gather}
where $t,q\in\{0,1,\ldots\},$ \ $p_1\cdc p_q\in\CC,$ \ and $p_q\ne0.$

\smallskip
We set
\begin{gather}
\label{h(la)} h(\la) =p_q^{-1}(\la^q+p_1\la^{q-1}+\ldots+p_q)
\end{gather}
and
\begin{gather}
\label{Z=h(A^k)} Z=h(A^u).
\end{gather}

\renewcommand{\labelenumi}{{\rm\theenumi)}}
\begin{proposition}
\label{p-Bprop} For any matrix $Z$ determined by \eqref{Z=h(A^k)}$,$ the following statements hold$:$
\begin{enumerate}
\item[\rm(1)] $ZC=CZ$ for any matrix $C$ commuting with $A;$
\item[\rm(2)] $Z$ is idempotent\/$;$
\item[\rm(3)] $AZ$ is nilpotent\/$;$
\item[\rm(4)] if $A$ is singular$,$ then $Z\ne0.$
\end{enumerate}
\end{proposition}

\propro{\ref{p-Bprop}}{ Item~(1) follows from the fact that $Z$ is a polynomial of~$A.$

(2) Let us consider the polynomial
\begin{gather*}
\varphi_0(\la)=\la h(\la)=p_q^{-1}\la(\la^q+p_1\la^{q-1}+\ldots+p_q).
\end{gather*}
Let $\widehat\psi(\la)$ be the minimal\x{} polynomial of\x{}~$A^u.$ Then, according to \cite[Theorem~3.3.1]{HornJohnson89}, $\varphi(\la)$ is divisible by $\widehat\psi(\la)$ and, by virtue of \cite[Theorem~3.3.6]{HornJohnson89} and Eq.\,\eqref{e-indA^k=1}, the lowest power\x{!!} of $\la$ in $\widehat\psi(\la)$ is at most~$\la^1.$\x{} Therefore, $\widehat\psi(\la)$ divides $\varphi_0(\la)$ and, consequently, $\varphi_0(\la)$ is an annihilating polynomial of\x{} $A^u,$ i.e., 
\begin{gather}
\label{A^kZ=0} 0=\varphi_0(A^u)=A^uh(A^u)=A^uZ=ZA^u.
\end{gather}
Let us prove that $Z^2=Z$. By definition, $h(\la)$ can be represented as $\la\,g(\la)+1,$ where $g(\la)$ is a polynomial in~\x{}$\la.$ Therefore, using\x{} \eqref{Z=h(A^k)} and \eqref{A^kZ=0} we obtain
\begin{gather*}
Z^2=Zh(A^u)=Z(A^ug(A^u)+I)=Z.\x{}
\end{gather*}

(3) If $u=0,$ then $Z=0$ by virtue of \eqref{A^kZ=0}, so in this case, $AZ$ is nilpotent with\x{} nilpotency index~$1.$ If $u>0,$ then by virtue of items~(1) and~(2) and\x{} Eq.\,\eqref{A^kZ=0} we obtain
\begin{gather*}
(AZ)^u=A^uZ^u=A^uZ=0,
\end{gather*}
and $AZ$ is nilpotent with nilpotency index\x{} at most~$u.$

(4) If $A$ is singular, then according to item~(c)\x{} of Theorem~4.6 in \cite{Ben-IsraelGreville74}, the lowest\x{} power of $\la$ in the minimal\x{} polynomial $\widehat\psi(\la)$ is $\la^{\ind A^u}\!,$\x{} consequently, it is not~$\la^0$. Therefore, $h(\la)$ is not divisible by $\widehat\psi(\la)$ and, consequently, $h(\la)$ is not an annihilating polynomial of\x{}~$A^\nu.$ Therefore,
$$
0\ne h(A^\nu)=Z.\eqno\qed
$$
}

\vspace{-1.2em}
The following theorem results from the 
condition~(f) of Section~\ref{s-chara} and Proposition~\ref{p-Bprop}.
\begin{theorem}
\label{t-eigen0} 
Let $\nu=\ind A$ and $u\ge\nu.$ Then for any nonzero annihilating polynomial of\x{} $A^u,$ $\varphi(\la),$ the matrix $Z$ defined by \eqref{Z=h(A^k)} is the eigenprojection of~$A.$
\end{theorem}

The eigenprojection of $A$ is related to\x{} the Drazin inverse $A^D$ (see (g) in Section~\ref{s-chara}). On the other hand,\x{} $A^D$~can be expressed\x{} via the eigenprojection $Z$ as follows%
\footnote[5]{For $\alpha=1,$ this expression was obtained in \cite{Rothblum76}.}
\cite{Koliha01}: 
\begin{gather}
\label{e-Z_Dr} A^D=(A+\alpha Z)^{-1}(I-Z),\quad\alpha\ne 0.
\end{gather}
Additionally, if $\,\ind A=1,$ then the group inverse $A^{\scriptscriptstyle\#} =A^D$ is representable as
\begin{gather*}
A^{\scriptscriptstyle\#}=(A+Z)^{-1}-Z
\end{gather*}
(see, e.g.,\x{} \cite{Meyer75,MeyerStadelmaier78}). Together with the formula $A^D=A^{\nu-1}{(A^\nu)}^{\scriptscriptstyle\#}$ which\x{} can be easily checked and the fact that $\ind A^\nu=1,$ the above representation allows one to obtain one more expression for~$A^D$:
\begin{gather*}
A^D=A^{\nu-1}((A^\nu+Z)^{-1}-Z).
\end{gather*}
It is worth mentioning\x{} that the Moore--Penrose generalized inverse\x{} can also be expressed in terms of the eigenprojection\x{} (see, e.g., \cite[pp.~686--687]{CharBen}).

\section{Determining the components and the minimal\x{} polynomial of~a~matrix}\label{s-co_mi}\vspace*{1mm}

Let $\la_1\cdc\la_s\in\CC$ be all distinct eigenvalues of~$A.$ We denote by $\nu_k$ the index of the eigenvalue $\la_k,$\; $k=1\cdc s,$ that is, the order of the greatest Jordan block corresponding to $\la_k,$ which is known to be equal to the power of the multiplier $(\la-\la_k)$ in the minimal\x{} polynomial of~$A.$ By the {\em eigenprojection of $A$ corresponding to $\la_k$} is meant the eigenprojection of the matrix $A-\la_kI.$

Having found\x{} $\la_1\cdc\la_s$ one can use Theorem~\ref{t-eigen0} to determine the corresponding eigenprojections. Let us\x{} turn to the theory of matrix components \cite{Gantmacher66,Lankaster78} according to which, for any function ${f:\CC\to\CC}$ having finite derivatives $f^{(j)}(\la_k)$\x{} of the corresponding orders $j=0\cdc\nu_k$ $(k=1\cdc s),$ the function $f(A)$ is naturally defined as follows:
\begin{gather*}
f(A):=\sum_{k=1}^s\sum_{j=0}^{\nu_k-1}f^{(j)}(\la_k)\,Z_{kj},
\end{gather*}
where under\x{} the derivative of order $0$ is meant the function itself\x{}.

The matrices $Z_{kj}$ are called the {\em components} of $A.$ Here, the component $Z_{k0}$ is the eigenprojection of $A$ corresponding to $\la_k,$\, $k=1\cdc s,$ the rest of the components being obtained by successive multiplication of $Z_{k0}$ by $A-\la_kI$ along with the factors\x{} $(j!)^{-1}$:
\begin{gather}
\label{e-next_comp} Z_{kj}=(j!)^{-1}(A-\la_kI)^j\,Z_{k0}
\end{gather}
(see \cite[Theorem~5.5.2]{Lankaster78}).

All components of a matrix are linearly independent (see, for example,\x{} \cite[\S~5.3]{Gantmacher66}), thus, in particular, there are no zero components among them. At the same time, $(A-\la_kI)^{\nu_k}\,Z_{k0}=0$ for each $k=1\cdc s,$ which follows from condition~(b) of Section~\ref{s-chara} and the fact that $Z_{k0}$ is the eigenprojection of $A-\la_kI,$ $\nu_k$ being its index. Therefore, after\x{} determining the component $Z_{k0}$ by means of Theorem~\ref{t-eigen0}, the remaining components can be found by multiplying $Z_{k0}$ by $A-\la_kI$ and the corresponding numerical factors\x{} (see~\eqref{e-next_comp}). As soon as the next multiplication by $A-\la_kI$ provides $0,$ all components corresponding to $\la_k$ are already determined and their number is~$\nu_k.$

In turn, knowledge of $\la_k$ and $\nu_k$ allows one to construct the minimal\x{} polynomial of $A$ as follows:
\begin{gather*}
\widehat\psi(\la)=\prod_{k=1}^s(\la-\la_k)^{\nu_k}.
\end{gather*}

\x{}Thus, starting with the eigenvalues of $A$ and\x{} annihilating polynomials for
\begin{gather*}
(A-\la_kI)^{u_k}\x{},\quad u_k\x{}\ge\ind(A-\la_kI),\quad k=1\cdc s,
\end{gather*}
one can proceed\x{} with determining the eigenprojections corresponding to all eigenvalues of $A,$ the components of $A,$ the indices of the eigenvalues, and the minimal\x{} polynomial of\x{}~$A.$

However, the following problem arises when this approach is implemented\x{} numerically. The above method involves verification of whether or not\x{} the current multiplication of $(A-\la_kI)^j\,Z_{k0}$ by $A-\la_kI$ results in~$0.$ Approximate calculations require\x{} an additional criterion for recognizing a matrix with ``small'' entries as a zero matrix. However, such\x{} problems are encountered\x{} in all algorithms for determining the rank of a matrix or the indices of its eigenvalues since\x{} these values are unstable to small perturbations of the matrix entries.

In the next section, we obtain explicit expressions for the eigenprojection and components of a square matrix whose eigenvalues are known.

\section{Eigenprojection and components of a matrix \\ Whose eigenvalues are known\x{}} \label{s-ep_eva}\vspace*{1mm}

One can readily check that if $J(\la)$ is a Jordan block corresponding to an eigenvalue $\la\ne0,$\x{} then $\ind(J(\la)-\la I) =\ind(J^p(\la)-\la^pI),$\, $p=2,3,\ldots$ (see, for example, \cite[Problem~16 in Section~3.2]{HornJohnson89}). This gives rise to the following lemma.
\vspace{-2mm}

\begin{lemma}
\label{l-ind} Let $\la_i\ne0$ be an\x{} eigenvalue of~$A.$ Then the index of $\la_i$ is equal to the index of $\la_i^p$ as the eigenvalue of $A^p,$\, $p=2,3,\ldots.$
\end{lemma}
\vspace{-2mm}

The following expressions\x{} for the eigenprojection and components of $A$ are\x{} based on knowing the eigenvalues of $A,$ rather than an annihilating polynomial of~$A^u,$ where $u\ge\ind A.$

\vspace{-2mm}
\begin{proposition}
\label{p-ZZk=} 
Let $A\in\CC^\nn;$ let $Z$ be the eigenprojection of~$A.$ Suppose that $\la_1\cdc\la_s$ are the distinct\x{} eigenvalues of~$A,$ $\nu_1\cdc\nu_s$ are their indices$,$ and the integers $u_1\cdc u_s$ are such that $u_i\ge\nu_i,$\, $i=1\cdc s.$ Let $u\ge\ind A.$ Then$:$
\begin{gather}
\label{e-Zeva} 
Z= \prod_{i:\,\la_i\ne0}\Bigl(I-(A/\la_i)^u\Bigr)^{\!u_i}
\end{gather}
\vspace{-1mm} and \vspace{-1mm}
\begin{gather}
\label{e-Zk_eva} Z_{kj}=
\prod_{i\ne k}\left(I-\left(\frac{A-\la_kI}{\la_i-\la_k}\right)^{\!u_k}\right)^{\!u_i}\!(j!)^{-1}(A-\la_kI)^j,
\end{gather}
where $Z_{kj}$ is the order $j$ component of $A$ corresponding to $\la_k,$\, $k=1\cdc s,$\, $j=0\cdc \nu_k-1.$
\end{proposition}
\vspace{-2mm}

If $j=0,$ then Eq.\,\eqref{e-Zk_eva} determines the eigenprojections of $A$ corresponding to its eigenvalues.

\propro{\ref{p-ZZk=}}{Note that
\begin{gather*}
\la\,\xi(\la)\equiv\la\prod_{i:\,\la_i\ne0}(\la-\la_i^u)^{u_i}
\end{gather*}
is an annihilating polynomial of\x{} $A^u.$ Indeed, it is divisible by the minimal\x{} polynomial of\x{}~$A^u$ because $\la_i^u$ are the eigenvalues of~$A^u,$ their indices do not exceed the numbers $u_i$ by Lemma~\ref{l-ind}, and $\ind A^u\le1$ according to \eqref{e-indA^k=1}. To prove \eqref{e-Zeva}, it suffices now to take $h(\la)$ (see \eqref{h(la)}) as the polynomial $\xi(\la)$ divided by its absolute term
\begin{gather*}
p=\prod_{i:\,\la_i\ne0}(-\la_i^u)^{u_i},\x{} 
\end{gather*}
apply Theorem~\ref{t-eigen0}, and perform some algebraic transformations.

By applying now \eqref{e-Zeva} to the matrix $A-\la_kI$ whose index, by definition, does not exceed $u_k,$ the distinct eigenvalues $\la_i-\la_k$ have the indices $\nu_i,$ $i=1\cdc s,$ respectively, and so $u_1\cdc u_s$ are non-strict upper bounds for these indices, we obtain \eqref{e-Zk_eva} with $j=0.$ Now the expressions \eqref{e-Zk_eva} for the components of $A$ of higher orders follow from \eqref{e-next_comp}.\qed}

Proposition~\ref{p-ZZk=} generalizes formula (5.4.3) in~\cite{Lankaster78} which refers to the case of $\nu_1=\ldots=\nu_s=1.$

\section*{Conclusion\x{}}\vspace*{1mm}

A number of necessary and sufficient conditions determining the eigenprojection (the principal idempotent) of a square matrix have been\x{} reviewed. A theorem which enables one to determine the eigenprojection of a matrix $A$ and thereby the Drazin inverse $A^D$ by means of any nonzero annihilating polynomial of $A^u$ with any $u\ge\ind A$ has been proved and used to determine the components of $A$ and its minimal\x{} polynomial. This theorem leads to explicit expressions for the eigenprojection and components of a matrix with known eigenvalues.\x{}

\revred{N.A.\ Bobylev}

\end{document}